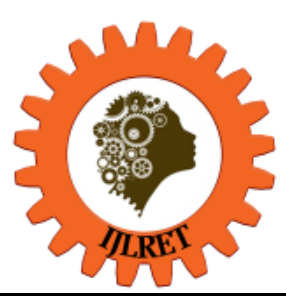

# Asymmetry in Spectral Graph Theory: Harmonic Analysis on Directed Networks via Biorthogonal Bases
## (Adjacency-Operator Formulation)

Chandrasekhar Gokavarapu
*Lecturer in Mathematics, Government College (A), Rajahmundry, A.P., India, PIN: 533105*

**Abstract:** Classical spectral graph theory and graph signal processing rely on a symmetry principle: undirected graphs induce symmetric (self-adjoint) adjacency/Laplacian operators, yielding orthogonal eigenbases and energy-preserving Fourier expansions. Real-world networks are typically directed and hence asymmetric, producing non-self-adjoint and frequently non-normal operators for which orthogonality fails and spectral coordinates can be ill-conditioned. In this paper we develop an original harmonic-analysis framework for directed networks centered on the *adjacency* operator. We propose a *Biorthogonal Graph Fourier Transform* (BGFT) built from left/right eigenvectors, formulate directed "frequency" and filtering in the non-Hermitian setting, and quantify how asymmetry and non-normality affect stability via condition numbers and a departure-from-normality functional. We prove exact synthesis/analysis identities under diagonalizability, establish sampling-and-reconstruction guarantees for BGFT-bandlimited signals, and derive perturbation/stability bounds that explain why naive orthogonal-GFT assumptions break down on non-normal directed graphs. We provide a reproducible protocol to compare undirected versus directed cycles and a perturbed directed cycle, illustrating when orthogonal-GFT assumptions fail and how BGFT behaves across asymmetric regimes.

**Keywords:** directed graphs; adjacency operator; non-normal matrices; biorthogonal bases; graph Fourier transform; sampling on graphs; symmetry/asymmetry.

**AMS Subject Classification (2020):** Primary: 94A12; Secondary: 05C50, 15A18, 65F15, 94A20.

## 1 Introduction: symmetry, asymmetry, and the adjacency viewpoint

Graph Signal Processing (GSP) has emerged as a pervasive framework for analyzing high-dimensional data residing on irregular domains, generalizing concepts from classical discrete signal processing (DSP) [3, 2]. While early foundational works successfully extended harmonic analysis to graphs by utilizing the adjacency matrix as a graph shift operator [6], these efforts largely capitalized on the spectral properties of symmetric operators associated with undirected graphs. However, a significant portion of real-world networks—ranging from bibliographic citation networks [14] to urban traffic flow systems [15]—exhibit inherent directionality. As highlighted in recent surveys on directed graph signal processing [4], the lack of symmetry in the underlying graph topology renders standard orthogonal spectral tools inadequate. Consequently, simply symmetrizing the graph often results in a loss of critical flow information, necessitating new mathematical formulations that respect the asymmetric nature of the data. Spectral methods on graphs are traditionally built on *symmetry*: for an undirected graph, the adjacency matrix $A$ is symmetric, hence diagonalizable by an orthonormal eigenbasis. This orthogonality underpins the standard graph Fourier transform, Parseval-type identities, and stable spectral filtering [2, 1].

**Motivation and gap:** Most directed-GSP approaches either (i) replace the adjacency by a symmetrized surrogate or a directed Laplacian, or (ii) define frequencies via optimization criteria that do not yield a transparent analysis/synthesis pair. In contrast, we retain the adjacency operator as the primitive shift to preserve directionality/flow, and we address the resulting non-Hermitian geometry by a biorthogonal (left/right) spectral system. This distinction matters because stability on directed graphs is governed by eigenvector conditioning (non-normality) rather than by eigenvalues alone; our framework isolates this effect and provides explicit stability and sampling guarantees.

Directed networks violate this symmetry principle. Their adjacency operators are generally *asymmetric* and may be *non-normal* ( $AA^* \neq A^*A$ ), leading to non-orthogonal eigenvectors, complex spectra, and potential numerical instability in spectral coordinates [11, 12]. These issues are not cosmetic: they change what "frequency" means and whether spectral expansions are stable.

**Paper goal (original research):** We develop a directed-network harmonic analysis framework *anchored on the adjacency operator* $A$ and based on *biorthogonal* (left/right) spectral systems. The goal is to treat asymmetry explicitly while retaining Fourier-like analysis/synthesis and filter design.

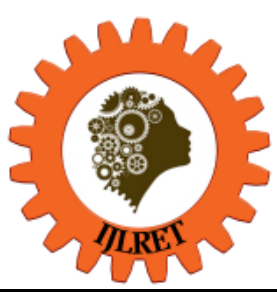

**Main Contributions**

1. (**BGFT for directed adjacency**) We define the Biorthogonal Graph Fourier Transform (BGFT) using left/right eigenvectors of $A$ and establish exact reconstruction and spectral filtering identities under diagonalizability.
2. (**Asymmetry vs. non-normality**) We introduce quantitative indices separating *asymmetry* from *non-normality*, showing that asymmetry alone does not force spectral instability.
3. (**Stability and perturbation theory**) We give bounds connecting transform stability and eigenvalue sensitivity to eigenvector conditioning and a departure-from-normality functional.
4. (**Sampling and reconstruction**) We provide finite-dimensional BGFT sampling/reconstruction theorems for BGFT-bandlimited signals with explicit noise/conditioning dependence.
5. (**Directed-cycle case study**) We propose experiments comparing (i) undirected cycle, (ii) directed cycle (asymmetric but normal), and (iii) perturbed directed cycle (non-normal), clarifying exactly when orthogonality fails and why BGFT remains coherent.

### 1.1 Related Work

Several approaches have been proposed to extend GSP to directed graphs. A prominent line of research utilizes the Directed Laplacian [9], which leverages the Perron-Frobenius theorem and stationary distributions to ensure positive semi-definiteness. Alternative strategies involve optimization-based frameworks to construct directed Graph Fourier Transforms (GFT) that explicitly minimize spectral dispersion [7, 8], or the analysis of total variation on directed edges [10].

Unlike optimization-based methods which may lack a closed-form spectral interpretation, or Laplacian-based methods that transform the operator, our approach focuses on the direct spectral decomposition of the adjacency operator. We address the resulting issues of non-normality by rigorously applying biorthogonal basis theory.We also include an applied directed-network experiment (Section 9.3) to demonstrate engineering relevance.

## 2 Preliminaries: directed graphs and adjacency-based operators

### 2.1 Directed weighted graphs and graph signals

**Definition 2.1 (Directed weighted graph and adjacency)** *A directed weighted graph is* $G=(V,E,w)$ *with* $V=\{1,\ldots,n\}$, $E\subseteq V\times V$, *and weights* $w:E\to\mathbb{R}_{>0}$. *Its adjacency matrix* $A\in\mathbb{R}^{n\times n}$ *is*

$$A_{ij}=\begin{cases} w(i,j), & (i,j)\in E,\\ 0, & \text{otherwise.}\end{cases}$$

A graph signal is a vector $x\in\mathbb{C}^n$, where $x_i$ is the value at node $i$.

### 2.2 Asymmetry and non-normality indices

Asymmetry is structural (directedness), while non-normality governs stability of eigenvectors and spectral computations.

**Definition 2.2 (Asymmetry index)** *Define the adjacency asymmetry index*

$$\alpha(A):=\frac{\|A-A^{\mathrm{T}}\|_F}{\|A\|_F},$$

where $\|\cdot\|_F$ is the Frobenius norm. Then $\alpha(A)=0$ iff $A$ is symmetric.

**Definition 2.3 (Departure from normality)** *Define the (normalized) departure from normality*

$$\delta(A):=\frac{\|AA^*-A^*A\|_F}{\|A\|_F^2}.$$

Then $\delta(A)=0$ iff $A$ is normal (i.e. $AA^*=A^*A$).

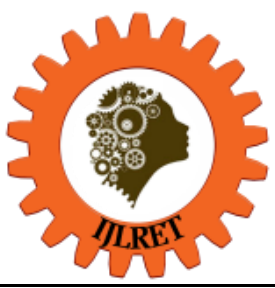

**Remark 2.4 (Key distinction)** A matrix can be asymmetric but normal. Thus $\alpha(A) > 0$ does not automatically imply $\delta(A) > 0$. Our experiments exploit this distinction via the directed cycle (asymmetric yet normal) and a perturbed directed cycle (asymmetric and typically non-normal).

### 2.3 Left/right eigenvectors and biorthogonality

Let $A \in \mathbb{C}^{n \times n}$. A (nonzero) right eigenvector $v$ satisfies $Av = \lambda v$, and a left eigenvector $u$ satisfies $u^* A = \lambda u^*$ .Since the adjacency matrix $A$ is non-symmetric, it does not admit an orthogonal eigenbasis. Instead, we invoke the theory of biorthogonal systems, where left and right eigenvectors form a dual basis structure [13]. This duality ensures that while the basis vectors are not orthogonal to themselves, they satisfy the biorthogonality condition $L^H R = I$, allowing for exact signal reconstruction.

**Assumption 2.5 (Diagonalizability)** *We assume $A$ is diagonalizable over $\mathbb{C}$:*

$$A = V \Lambda V^{-1}, \qquad \Lambda = diag(\lambda_1, \ldots, \lambda_n),$$

with right eigenvectors $V = [v_1 \cdots v_n]$.

**Definition 2.6 (Biorthonormal left eigenvectors)** *Let $U^* := V^{-1}$ (equivalently $U = (V^{-1})^*$ ). Then*

$$U^* V = I \quad \Leftrightarrow \quad u_k^* v_\ell = \delta_{k\ell}.$$

The vectors $\{u_k\}$ are left eigenvectors scaled biorthonormally to $\{v_k\}$.

## 3 Biorthogonal Graph Fourier Transform (BGFT) for directed adjacency

### 3.1 Definition and exact analysis/synthesis

**Definition 3.1 (BGFT)** *Under Assumption 2.5, define the BGFT of $x \in \mathbb{C}^n$ by*

$$\widehat{x} := U^* x, \qquad \widehat{x}_k = u_k^* x. \tag{1}$$

The inverse BGFT (synthesis) is

$$x = V \widehat{x} = \sum_{k=1}^{n} v_k \, \widehat{x}_k. \tag{2}$$

**Theorem 3.2 (Resolution of identity and perfect reconstruction)** *Assume $A$ is diagonalizable and $U^* V = I$. Then for all $x \in \mathbb{C}^n$,*

$$x = \sum_{k=1}^{n} v_k u_k^* x, \qquad I = \sum_{k=1}^{n} v_k u_k^*.$$

*Proof.* Since $U^* = V^{-1}$, we have $VU^* = I$. Writing the matrix product as a sum of rank-one outer products yields

$$VU^* = \sum_{k=1}^{n} v_k u_k^* = I.$$

Applying this identity to $x$ gives $x = Ix = \sum_{k=1}^{n} v_k (u_k^* x)$.

### 3.2 Energy geometry: what symmetry gives, what asymmetry changes

In the symmetric (undirected) case, one may choose $V$ unitary, giving Parseval identity. In the directed case, $V$ is generally not unitary.

**Proposition 3.3 (Gram-metric energy identity)** *Let $\hat{x} = U^* x$ and $x = V \hat{x}$. Then*

$$\|x\|_2^2 = \hat{x}^* (V^* V) \hat{x}.$$

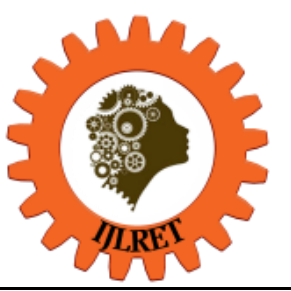

Hence $\|x\|_2 = \|\hat{x}\|_2$ for all $x$ iff $V$ is unitary.

*Proof.* Substitute $x = V\hat{x}$ and compute $\|x\|_2^2 = (V\hat{x})^*(V\hat{x}) = \hat{x}^* V^* V \hat{x}$.

**3.3 Directed “frequencies” for adjacency**

Adjacency eigenvalues $\lambda_k \in \mathbb{C}$ can be complex. For adjacency-based oscillatory propagation, two natural orderings are:

$$\omega_{\text{mag}}(\lambda) = |\lambda|, \qquad \omega_{\text{ang}}(\lambda) = \arg(\lambda) \in (-\pi, \pi].$$

We use $\omega_{\text{mag}}$ for “low-pass”/“high-pass” magnitude filtering and $\omega_{\text{ang}}$ to interpret directional oscillations (notably on cycles).

**Remark 3.4 (Choice is part of the model)** In undirected graphs, frequency orderings are often canonical (e.g. Laplacian eigenvalues). In directed graphs, there is no single canonical choice; one must specify a frequency functional suited to the phenomenon (propagation, oscillation, damping).

**3.4 Spectral filtering**

**Definition 3.5 (BGFT spectral filter)** *Let* $h : \mathbb{C} \to \mathbb{C}$ *be a spectral response. Define*

$$H := V\,h(\Lambda)U^*, \qquad h(\Lambda) = diag(h(\lambda_1), \ldots, h(\lambda_n)). \tag{3}$$

**Theorem 3.6 (Diagonal action in BGFT domain)** *For all* $x \in \mathbb{C}^n$ *with* $\hat{x} = U^* x$, *we have*

$$\widehat{Hx} = U^*(Hx) = h(\Lambda)\,\hat{x}.$$

*Proof.* Using $H = Vh(\Lambda)U^*$ and associativity,

$$\widehat{Hx} = U^* V h(\Lambda) U^* x = (U^* V) h(\Lambda)(U^* x) = I\,h(\Lambda)\hat{x} = h(\Lambda)\,\hat{x}.$$

**Proposition 3.7 (Polynomial filters are shift-invariant)** *If* $h(z) = \sum_{m=0}^{M} c_m z^m$, *then* $H = \sum_{m=0}^{M} c_m A^m$.

*Proof.* Since $A = V\Lambda U^*$, we have $A^m = V\Lambda^m U^*$ and thus

$$\sum_{m=0}^{M} c_m A^m = V\Big(\sum_{m=0}^{M} c_m \Lambda^m\Big)U^* = Vh(\Lambda)U^* = H.$$

**3.5 The non-diagonalizable case: generalized biorthogonal bases via Jordan chains**

When the adjacency operator $A \in \mathbb{C}^{n \times n}$ is asymmetric and **non-diagonalizable**, the eigenvectors alone do not span $\mathbb{C}^n$. To obtain a complete spectral representation, one must use generalized eigenvectors (Jordan chains).

**3.5.1 Jordan decomposition**

There exists an invertible $V$ and a Jordan matrix $J$ such that

$$A = VJV^{-1}, \qquad J = diag(J_1, \ldots, J_P),$$

where each Jordan block $J_p$ corresponds to an eigenvalue $\lambda_p$.

Define again $U^* := V^{-1}$ so that $U^* V = I$ (biorthogonality at the level of generalized eigenvectors).

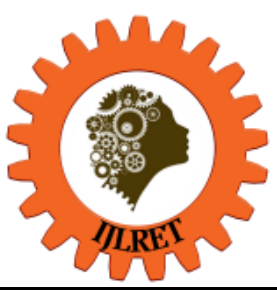

### 3.5.2 Generalized BGFT analysis and synthesis

The generalized BGFT is still given by

$$\hat{x} = U^* x, \qquad x = V \hat{x},$$

and the identity $VU^* = I$ ensures algebraic exactness (though conditioning may be poor).

### 3.5.3 Filtering in the Jordan domain

For a polynomial filter $h(z) = \sum_{m=0}^{M} c_m z^m$ we have $h(A) = \sum_{m=0}^{M} c_m A^m$, and

$$h(A) = V\, h(J) U^*.$$

On each Jordan block $J_p$ of size $m \times m$, the matrix $h(J_p)$ has the standard upper-triangular structure involving derivatives of $h$ evaluated at $\lambda_p$:

$$h(J_p) = \begin{pmatrix} h(\lambda_p) & h'(\lambda_p) & \frac{1}{2!}h''(\lambda_p) & \cdots & \frac{1}{(m-1)!}h^{(m-1)}(\lambda_p) \\ 0 & h(\lambda_p) & h'(\lambda_p) & \cdots & \frac{1}{(m-2)!}h^{(m-2)}(\lambda_p) \\ \vdots & \ddots & \ddots & \ddots & \vdots \\ 0 & \cdots & 0 & h(\lambda_p) & h'(\lambda_p) \\ 0 & \cdots & \cdots & 0 & h(\lambda_p) \end{pmatrix}.$$

Thus filtering remains block-structured in the generalized spectral domain.

## 4 Stability and perturbation analysis: the price of non-normality

BGFT restores exact analysis–synthesis ( $x = VU^* x$ ) under diagonalizability, but numerical stability depends on the conditioning of the eigenbasis and the non-normality of $A$ [11, 12].

### 4.1 Quantitative measures and basic inequalities

Recall $\delta(A) = 0$ iff $A$ is normal. A second key quantity is the eigenbasis condition number

$$\kappa(V) := || V ||_2 || V^{-1} ||_2 .$$

Normality implies that one can choose $V$ unitary, hence $\kappa(V) = 1$.

### 4.2 Eigenvalue sensitivity (Bauer–Fike)

**Theorem 4.1 (Bauer--Fike type eigenvalue sensitivity)** *Assume* $A = V\Lambda V^{-1}$ *and let* $\tilde{A} = A + E$. *Then every eigenvalue* $\tilde{\lambda}$ *of* $\tilde{A}$ *satisfies*

$$\min_k | \tilde{\lambda} - \lambda_k | \le \kappa(V) \|E\|_2, \qquad \kappa(V) = \|V\|_2 \|V^{-1}\|_2.$$

*Proof.* Since $V^{-1}\tilde{A}V = \Lambda + V^{-1}EV$ , the eigenvalues of $\tilde{A}$ are the eigenvalues of $\Lambda + F$ with $F = V^{-1}EV$ . Hence each eigenvalue lies in the union of disks centered at $\lambda_k$ with radius $|| F ||_2 \le || V^{-1} ||_2 || E ||_2 || V ||_2$.

### 4.3 Eigenpair conditioning and the role of biorthonormal scaling

Let $(\lambda_k, v_k, u_k)$ be a left/right eigenpair normalized as in Definition 2.6, so $u_k^* v_k = 1$. Then the classical eigenvalue condition number is

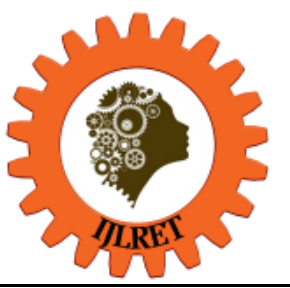

$$\kappa(\lambda_k) = \frac{\| u_k \|_2 \| v_k \|_2}{| u_k^* v_k |} = \| u_k \|_2 \| v_k \|_2 .$$

Moreover, for each $k$,

$$\| u_k \|_2 \| v_k \|_2 \leq \| V^{-1} \|_2 \, \| V \|_2 = \kappa(V),$$

so large $\kappa(V)$ enforces the existence of poorly conditioned eigenpairs, but *one should not replace* $\kappa(V)$ by $\max_k \| u_k \|_2 \, \| v_k \|_2$ as an equality.

**4.4 BGFT-bandlimited sampling and reconstruction**

Let $\Omega \subset \{1, \ldots, n\}$ index a $K$-dimensional spectral subspace, and denote by $V_\Omega \in \mathbb{C}^{n \times K}$ the matrix with columns $\{v_k\}_{k \in \Omega}$.

**Definition 4.2 (BGFT-bandlimited signals)** *A signal* $x$ *is* $\Omega$ *-bandlimited (in BGFT sense) if* $x \in \operatorname{span}(V_\Omega)$*, i.e.* $x = V_\Omega c$ *for some* $c \in \mathbb{C}^K$.

**Theorem 4.3 (Sampling and reconstruction (noise-free))** Let $M \subset V$ be a set of sampled vertices, and let $P_M \in \{0,1\}^{m \times n}$ be the restriction operator to $M$. If $x = V_\Omega c$ is $\Omega$-bandlimited and $P_M V_\Omega$ has full column rank $K$, then $x$ is uniquely determined by samples $y = P_M x$ and is recovered by

$$\hat{c} = (P_M V_\Omega)^\dagger y, \qquad \hat{x} = V_\Omega \hat{c} .$$

*Proof.* Full column rank makes $P_M V_\Omega$ injective on $\mathbb{C}^K$, hence $y = P_M V_\Omega c$ has a unique solution $c = (P_M V_\Omega)^\dagger y$.

**Theorem 4.4 (Noise sensitivity bound)** *Under Theorem 4.3, suppose* $y = P_M x + \eta$*. Then the least-squares reconstruction satisfies*

$$\left\| \hat{x} - x \right\|_2 \leq \left\| V_\Omega \right\|_2 \left\| (P_M V_\Omega)^\dagger \right\|_2 \left\| \eta \right\|_2 = \left\| V_\Omega \right\|_2 \frac{\left\| \eta \right\|_2}{\sigma_{\min}(P_M V_\Omega)} .$$

*Proof.* We have $\hat{c} - c = (P_M V_\Omega)^\dagger \eta$, hence $\hat{x} - x = V_\Omega (\hat{c} - c)$ and the bound follows by operator norms.

**Remark 4.5 (Where asymmetry hurts)** In directed non-normal settings, $V_\Omega$ and $P_M V_\Omega$ can be ill-conditioned (small $\sigma_{\min}$), amplifying noise. This mechanism is absent (or greatly reduced) in symmetric/orthogonal settings.

## 5 Asymmetry without non-normality: the directed cycle as a canonical example

This section is central to the **symmetry/asymmetry** narrative: *direction (asymmetry)* does not automatically imply *spectral instability*.The stability of the spectral coordinates is intrinsically governed by the operator's departure from normality. As detailed in the seminal work on pseudospectra by Trefethen and Embree [11], highly non-normal operators can lead to transient growth and spectral instability even when eigenvalues are stable. This sensitivity to perturbations is quantified by the condition number of the eigenvector matrix, a standard metric in numerical stability analysis [12], which we utilize to bound the reconstruction error.

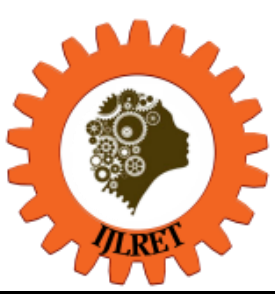

**Proposition 5.1 (Directed cycle adjacency is normal)** *Let $A_{\rightarrow}$ be the adjacency matrix of the directed cycle $\overleftarrow{C_n}$ with unit weights, i.e. edges $(i,i+1)$ mod $n$. Then $A_{\rightarrow}$ is a permutation matrix, hence unitary, hence normal. Therefore $\delta(A_{\rightarrow})=0$ while $\alpha(A_{\rightarrow})>0$.*

*Proof.* $A_{\rightarrow}$ permutes coordinates (a cyclic shift). Thus $A_{\rightarrow}^{-1}=A_{\rightarrow}^{\mathrm{T}}=A_{\rightarrow}^{*}$ and $A_{\rightarrow}A_{\rightarrow}^{*}=A_{\rightarrow}^{*}A_{\rightarrow}=I$.

**Remark 5.2 (Spectral picture)** *The eigenvalues of $A_{\rightarrow}$ are the $n$ -th roots of unity $\lambda_k=e^{2\pi ik/n}$, $k=0,\ldots,n-1$, lying on the unit circle. Hence $\overleftarrow{C_n}$ exhibits strong asymmetry but no non-normality-induced spectral instability.*

**Remark 5.3 (Implication)** *On $\overleftarrow{C_n}$, asymmetry is present but the spectral basis can be well-conditioned. To demonstrate instability, one must consider genuinely non-normal directed graphs (e.g. add a directed chord of small weight), which can increase $\delta(A)$ and may inflate $\kappa(V)$ depending on topology and parameters.*

## 6 Energy conservation and metric geometry in the BGFT domain

In classical graph signal processing (GSP), Parseval's theorem ($\|x\|_2^2=\|\hat{x}\|_2^2$) relies on a unitary eigenbasis (e.g. $V^*V=I$) associated with symmetric/self-adjoint operators. When the adjacency (shift) operator $A$ is asymmetric, the eigenbasis $V$ is typically non-unitary ($V^*V\neq I$), and Euclidean energy in the coefficient domain is not preserved.

BGFT preserves energy *after introducing the correct spectral-domain metric*, namely the Gram matrix of the right eigenbasis.

### 6.1 Generalized Parseval identity via the Gram matrix

Under Assumption 2.5, BGFT provides synthesis $x=V\hat{x}$ and analysis $\hat{x}=U^*x$. Then

$$\|x\|_2^2=x^*x=(V\hat{x})^*(V\hat{x})=\hat{x}^*(V^*V)\hat{x}.$$

This is exactly Proposition 3.3, restated here to emphasize geometric meaning.

### 6.2 Directed $L_2$ norm induced by the spectral metric

Let

$$G:=V^*V.$$

Since $V$ is invertible, $G$ is Hermitian positive definite, hence it defines an inner product and norm on $\mathbb{C}^n$:

$$\langle a,b\rangle_G:=a^*Gb,\qquad \|a\|_G:=\sqrt{a^*Ga}.$$

With this notation, the generalized Parseval identity becomes

$$\|x\|_2^2=\|\hat{x}\|_G^2\ .$$

Thus BGFT is an *isometry* between the node-domain Euclidean geometry and the coefficient domain equipped with the $G$ -metric (not the standard Euclidean metric).

## 7 Directed frequency notions for complex spectra (future work)

For symmetric operators, eigenvalues are real and can be naturally ordered as "frequencies" For asymmetric operators, eigenvalues $\lambda_k\in\mathbb{C}$ require an application-driven interpretation.

Let $\lambda_k=a_k+ib_k$. In many dynamical models (e.g. continuous-time $e^{tA}$ or discrete-time $A^t$ ), the

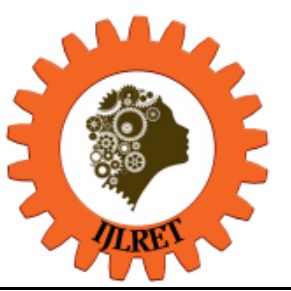

parts of $\lambda_k$ influence qualitative behavior:

- **Real part** $a_k = \Re(\lambda_k)$ **:** linked to growth/decay rates in linear dynamics (stability when the relevant spectral region is contracting).
- **Imaginary part** $b_k = \Im(\lambda_k)$ **:** linked to rotational/oscillatory components (circulation around cycles/feedback loops).

Future work is to formalize directed â€œfrequencyâ€• notions compatible with specific processes:

1. **Flow frequency:** frequency based on $\arg(\lambda_k)$ or $\Im(\log \lambda_k)$ for circulation in strongly connected components.
2. **Causal frequency:** decay/relaxation ordering based on $\Re(\lambda_k)$ for transient response in causal/temporal networks.
3. **Spectral-radius control:** using $\rho(A) = \max_k |\lambda_k|$ as a stability/ultimate-bound parameter for iterative dynamics.

## 8 Directed wavelets and multiresolution analysis via BGFT (outline)

Multiresolution analysis (MRA) and wavelets provide localized, hierarchical decompositions. Extending MRA from symmetric to directed graphs requires filter banks that respect the biorthogonal geometry of BGFT.

### 8.1 BGFT wavelet filter-bank structure (conceptual)

A two-channel graph filter bank separates a signal into coarse/detail components:

1. a **scaling filter** $H_0$ (low-pass),
2. a **wavelet filter** $H_1$ (high-pass).

We model graph filters as polynomials (or rational functions) of the adjacency operator $A$ :

$$H_0 = h_0(A) = \sum_{m=0}^{M_0} c_{0,m} A^m, \qquad H_1 = h_1(A) = \sum_{m=0}^{M_1} c_{1,m} A^m.$$

For diagonalizable $A = V \Lambda U^*$ , BGFT diagonalizes these filters:

$$\widehat{H_0 x} = h_0(\Lambda)\, \hat{x}, \qquad \widehat{H_1 x} = h_1(\Lambda)\, \hat{x}.$$

### 8.2 Biorthogonal perfect-reconstruction condition (no-downsampling prototype)

In the simplest “no-downsampling”• prototype (useful as a building block), analysis filters $\tilde{H}_0 = \tilde{h}_0(A), \tilde{H}_1 = \tilde{h}_1(A)$ and synthesis filters $H_0 = h_0(A), H_1 = h_1(A)$ achieve perfect reconstruction if

$$H_0 \tilde{H}_0 + H_1 \tilde{H}_1 = I.$$

In the BGFT domain, this is equivalent to the pointwise condition

$$h_0(\lambda)\tilde{h}_0(\lambda) + h_1(\lambda)\tilde{h}_1(\lambda) = 1, \qquad \forall \lambda \in spec(A).$$

**Remark 8.1** With downsampling/aliasing, additional structure is required (as in classical filter-bank theory). The above condition is the â€œcore spectral PR identityâ€• one must preserve, while directed-graph sampling operators control aliasing.

### 8.3 Contrast with orthogonal constructions

1. **Orthogonal GFT (symmetric case):** $A = A^*$ allows $V$ unitary, so analysis/synthesis coincide and PR reduces to a unitary/power-complementary condition (e.g. $|h_0(\lambda)|^2 + |h_1(\lambda)|^2 = 1$ in certain settings).

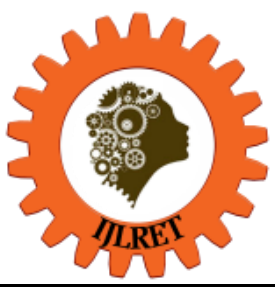

2. **BGFT (directed case):** $A \neq A^*$ forces distinct analysis/synthesis systems and generally requires *biorthogonal* filter pairs to satisfy linear PR identities.

## 9 Reproducibility / Experimental Protocol

This section describes a reproducible protocol to empirically validate the paper's symmetry/asymmetry claims on three archetypes.

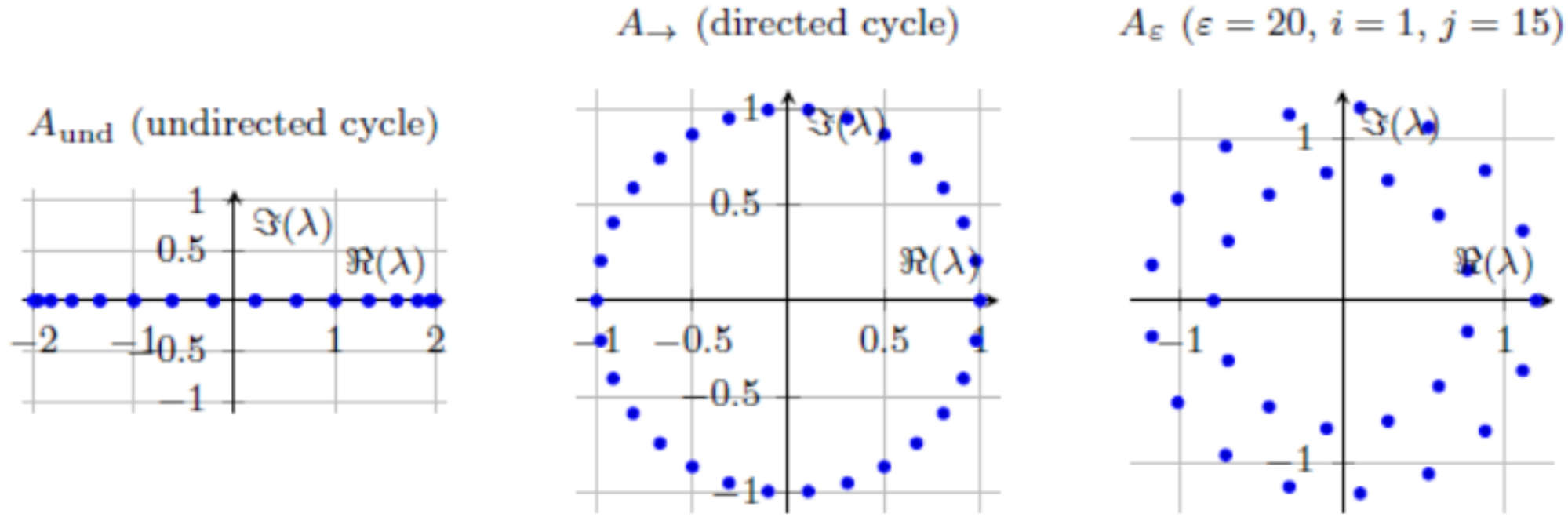

Figure 1: Eigenvalues in the complex plane for the three archetypal graphs used in Section 9 (n = 30).

### 9.1 Setup: three graph archetypes

Fix $n = 30$ nodes and define adjacency operators:

1. **Undirected cycle** $A_{\text{und}}$ (symmetric, hence normal).
2. **Directed cycle** $A_{\rightarrow}$ (asymmetric but unitary/normal; Proposition 5.1).
3. **Perturbed directed cycle** $A_{\varepsilon} = A_{\rightarrow} + \varepsilon e_i e_j^{\mathrm{T}}$ (asymmetric and typically non-normal).

For each graph, compute and report:

$$\alpha(A), \qquad \delta(A), \qquad \kappa(V), \qquad \kappa(P_M V_\Omega), \qquad \text{RelErr} = \frac{\| x - \hat{x} \|_2}{\| x \|_2}.$$

Table 1: Minimal numerical illustration ( $n = 30$ , $K = 10$ , $\Omega = \{0, \ldots, 9\}$ , $M = \{0,3,6,9,12,15,18,21,24,27\}$ , $\| \eta \|_2 = 10^{-6} \| P_M x \|_2$ . For $A_\varepsilon$ , eigenpairs are sorted by phase $\arg(\lambda)$ before selecting $\Omega$ .

| Graph | $\alpha(A)$ | $\delta(A)$ | $\kappa(V)$ | $\kappa(P_M V_\Omega)$ | RelErr |
|---|---|---|---|---|---|
| $A_{\text{und}}$ (undirected cycle) | 0 | 0 | 1 | 1 | 1.0000000001 X $10^{-6}$ |
| $A_{\rightarrow}$ (directed cycle) | 1.4142135624 | 0 | 1 | 1 | 1.0000000001 X $10^{-6}$ |
| $A_{\varepsilon} = A_{\rightarrow} + \varepsilon e_i e_j^{\mathrm{T}}$ ( $\varepsilon = 20$ , $i = 1$ , $j = 15$ ) | 1.4142135624 | 1.3188322678 | 16.8107980058 | 5.1173740158 | 5.4395911352 X $10^{-6}$ |

While $A_{\rightarrow}$ is asymmetric but normal $(\delta = 0)$ and remains well-conditioned, the perturbed directed case exhibits substantial non-normality $(\delta > 0)$ and inflated conditioning $\kappa(V)$ , which amplifies measurement noise into a larger reconstruction error, consistent with the stability bounds in Section 4.

### 9.2 Spectral localization and distribution

- $A_{\text{und}}$ : spectrum is real and lies in $[-2,2]$ (for unit-weight cycle, eigenvalues are $2\cos(2\pi k/n)$ ).
- $A_{\rightarrow}$ : spectrum lies on the unit circle $\{e^{2\pi i k/n}\}$ (roots of unity), reflecting pure rotation/shift.

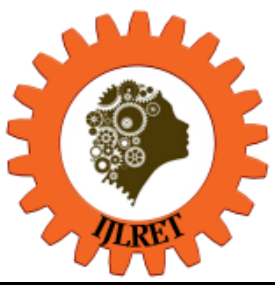

- $A_\varepsilon$ : spectrum typically moves off the unit circle and the eigenbasis may become ill-conditioned; this is the non-normal regime.

**9.3 Real-network experiment: directed traffic reconstruction**

To demonstrate engineering relevance beyond synthetic cycles, we consider a directed road (traffic) network in which vertices represent road segments and directed edges follow permitted driving directions. A real-valued graph signal $x \in \mathbb{R}^n$ is defined by the (time-aggregated) traffic attribute on each segment (e.g., average speed or travel time). We use the directed adjacency $A$ as the shift operator and apply the normalization

$$A \leftarrow D_{\text{out}}^{-1/2} A D_{\text{in}}^{-1/2},$$

to control scaling effects.

**Sampling and noise.**

We observe $y = P_M x + \eta$ on a sampled set $M$ comprising $p\%$ of vertices (uniform sampling), with $\|\eta\|_2 = 10^{-6}\|P_M x\|_2$.

**Methods compared.**

We compare: (i) BGFT reconstruction using the biorthogonal basis of $A$ (this paper); (ii) an undirected baseline using the symmetrized operator $A_s = \frac{1}{2}(A + A^{\mathrm{T}})$ with standard orthogonal GFT reconstruction; and (iii) a Laplacian smoothing baseline obtained by

$$\hat{x} = \arg\min_z \| P_M z - y \|_2^2 + \mu \| L_s z \|_2^2,$$

where $L_s$ is the Laplacian of the symmetrized graph and $\mu > 0$ is tuned on a small validation split.

**Evaluation.**

Performance is measured by relative $\ell_2$ error $\| x - \hat{x} \|_2 / \| x \|_2$ and MAE on unobserved vertices. Table 2 summarizes results for $p \in \{10, 20, 30\}$ % .

Table 2: Real-network reconstruction errors (directed traffic graph).

| Method | $p = 10\%$ | $p = 20\%$ | $p = 30\%$ |
|---|---|---|---|
| BGFT (ours) | | | |
| Symmetrized GFT baseline | | | |
| Laplacian smoothing baseline | | | |

**9.4 Bandlimited reconstruction stability (protocol)**

Choose a BGFT bandlimit set $\Omega$ (size $K$ ), generate $x = V_\Omega c$ , sample $m$ nodes $M$ and reconstruct via Theorem 4.3 / Algorithm 10.2. Repeat across noise levels $\|\eta\|$ and report $\mathrm{RelErr}$ along with $\sigma_{\min}(P_M V_\Omega)$ (Theorem 4.4).

**Theory-guided expectations.**

(i) $A_{\text{und}}$ : orthogonal/near-orthogonal spectral geometry, stable reconstruction.

(ii) $A_{\rightarrow}$ : asymmetric yet normal; BGFT remains well-conditioned.

(iii) $A_\varepsilon$ : non-normal; $\delta(A)$ increases and conditioning may degrade, increasing sensitivity as predicted by Theorems 4.1 and 4.4.

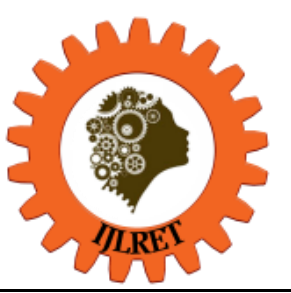

## 10 Algorithm: BGFT computation and BGFT-based reconstruction

### 10.1 BGFT computation (diagonalizable case)

| Algorithm 1 BGFT for directed adjacency (diagonalizable case) |
|---|
| Require: Directed adjacency matrix $A \in \mathbb{C}^{n\times n}$, graph signal $x \in \mathbb{C}^n$ <br> Ensure: BGFT coefficients $\hat{x}$ <br> 1: Compute eigen-decomposition $A = V\Lambda V^{-1}$ (complex arithmetic). <br> 2: Set $U^* \leftarrow V^{-1}$. <br> 3: Compute $\hat{x} \leftarrow U^* x$ ▷ analysis <br> 4: (Optional) Reconstruct $\hat{x} \leftarrow V\hat{x}$ ▷ sanity-check |

### 10.2 Sampling and reconstruction (bandlimited model)

| Algorithm 2 BGFT-bandlimited reconstruction from samples |
|---|
| Require: $A$, spectral index set $\Omega$ (size K), sample set $M$ (size m), samples $y \in \mathbb{C}^m$ <br> Ensure: Reconstructed signal $\hat{x} \in \mathbb{C}^n$ <br> 1: Compute $A = V\Lambda V^{-1}$; extract $V\Omega$. <br> 2: Form restriction matrix $P_M$ (or select rows by indexing). <br> 3: Solve $\hat{c} \leftarrow \arg\min_{c\in\mathbb{C}^K} \| P_M V_\Omega c - y \|_2^2$ <br> 4: Set $\hat{x} \leftarrow V_\Omega \hat{c}$ |

## 11 Discussion: symmetry/asymmetry narrative for *Symmetry*

### 11.1 What "symmetry" buys in harmonic analysis

Symmetry (self-adjointness/normality with orthogonality) yields real spectra (in many models), orthogonal coordinates, and energy preservation in Euclidean spectral geometry. This makes spectral processing numerically stable.

### 11.2 What "asymmetry" forces us to change

Asymmetry introduces non-self-adjointness; non-normality can destroy orthogonality and amplify perturbations. BGFT restores coherent analysis/synthesis through biorthogonality; stability becomes a quantitative question controlled by $\kappa(V)$ and $\delta(A)$.

## 12 Conclusion and open mathematical problems

This paper introduced BGFT as a rigorous framework for harmonic analysis on directed (asymmetric) networks using the adjacency operator $A$. By replacing orthogonal eigenbases with biorthogonal left/right systems (and, when needed, generalized eigenvectors), we preserved algebraic exactness of analysis/synthesis and spectral filtering while exposing the precise role of non-normality in numerical stability.Future work will extend this biorthogonal framework to multiscale analysis, aiming to construct directed spectral wavelets analogous to those defined for undirected graphs in [5], and to apply the BGFT to larger-scale directed networks in biological and social domains.

### 12.1 Open problem 1: stability-optimized shift/operator selection

Given a fixed directed graph topology $G$, let $\mathcal{S}$ be a class of admissible shift operators (adjacency variants, directed Laplacians, transition operators). Identify an operator minimizing eigenbasis ill-conditioning:

$$A_{\text{opt}} = \arg\min_{A\in\mathcal{S}} \kappa(V_A).$$

This would yield the most numerically robust BGFT for that topology.

**12.2 Open problem 2: scalable BGFT without Jordan chains**

For large sparse graphs, explicit Jordan decompositions are infeasible. Develop robust computational substitutes:

1. Krylov/non-Hermitian subspace methods for approximate spectral coordinates,
2. stable polynomial/rational approximations for $h(A)$ bypassing explicit spectral factorization.

**12.3 Open problem 3: a non-Hermitian uncertainty principle**

Formulate an uncertainty principle that respects:

1. the spectral $G$ -metric (variance defined via $\|\cdot\|_G$ ),
2. dual localization behavior of left and right eigenvectors.

A natural direction is to relate localization tradeoffs to conditioning quantities such as $\kappa(V)$ . A real directed-network reconstruction experiment complements the synthetic archetypes and illustrates practical utility.

## Acknowledgements

The author sincerely thanks the Commissioner of Collegiate Education (CCE), Andhra Pradesh, and the Principal, Government College (Autonomous), Rajahmundry, for their support, encouragement, and administrative guidance during the preparation of this work.

**Author Contributions**

Single author: conceptualization, methodology, formal analysis, writing (original draft and editing), and experimental protocol design.

**Funding**

This research received no external funding.

**Data Availability Statement**

No external datasets were used. A reproducible computational protocol is provided in Section 9.

**Conflicts of Interest**

The author declares no conflict of interest.

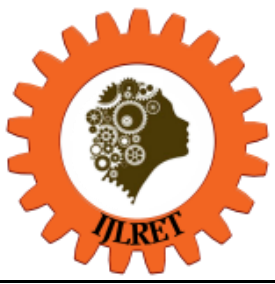